\magnification=\magstep2  
\hsize=16 truecm
\vsize=24 truecm
\baselineskip = 8 truemm  
\parindent=1.2truecm      
\overfullrule=0pt  
\language=0        
\footline={\hfil\ifnum\pageno=1\else\rm\the\pageno\fi\hfil}
\font\BF=cmbx12  
\font\cmbsy=cmbsy10
\def\Par{\hbox{\cmbsy\char120}}  
\font\rsfs=rsfs10
\def\a{\alpha}  \def\b{\beta}  \def\g{\gamma}  \def\G{\Gamma}  \def\t{\theta}
\def\s{\sigma} \def\u{\tau}
\def\mod{\mathop{\hbox{\rm mod}}}
\def\reg{\mathop{\hbox{\rm reg}}}
\def\Re{\mathop{\hbox{\rm Re}}}
\font\msbm=msbm10
\def\RR{\hbox{\msbm R}}  
\def\CC{\hbox{\msbm C}}  
\def\QQ{\hbox{\msbm Q}}  
\def\KK{\hbox{\msbm K}}  
\def\ZZ{\hbox{\msbm Z}}
\font\cyr=cmcyr10
\def\No{\hbox{\cyr\char'031\hskip2pt}}   

\def\sgn{\mathop{\hbox{\rm sgn}}}
\def\Sp{\mathop{\hbox{\rm Sp}}}
\def\script#1{\hbox{\rsfs #1}}

\vskip 12 truemm
\centerline{\BF Adelic Formulas for Gamma- and Beta-functions}
\centerline{\BF in Algebraic Numbers Fields}
\vskip 9 truemm
\centerline{\BF V.~S.~Vladimirov}
\medskip
\centerline {(Steklov Mathematical Institute, Moscow, Russia)}
\bigskip
{\bf Abstract.} On the basis of analysis on the adele ring of any algebraic
numbers field (Tate's formula) a regularization for divergent adelic products 
of gamma- and beta-functions for all completions of this field are proposed, 
and cor-ponding regularized adelic formulas are obtained.
\bigskip
\centerline {\BF \Par 1. Introduction}
\medskip

Freund and Witten~[1](1987) and, independently a little later, Volovich~
[2] (1988), suggested for the beta-function $B_\infty$ of the real numbers 
field $\RR$ the following adelic formula
$$B_\infty (\a,\b)\prod_{p=2}^\infty B_p(\a,\b)=1. \eqno (1.1)$$
where $B_p$ is the beta-function of p-adic numbers field $\QQ_p$.

However, the infinite product in (1.1) diverges for all $\a$ and $\b$ and,
therefore, it is not clear how to understand Eq.~(1.1). A detailed
discussions of these matters can be found in~[3]--[8].

In this paper on the basis of analysis on the adele ring (Tate's formular)
for any algebraic numbers field a regularization for divergent infinite 
products of gamma- and beta-functions for all characters (ramified or non) is 
proposed, and the regularized adelic formulas are derived.

More precisely, let $\KK=\QQ(\epsilon)$ be an algebraic numbers field of degree
$n=\s+2\u$ where $\s$ and $2\u$ are numbers of real and complex roots resp. of
the minimal polynomial for a primitive element $\epsilon$ of field $\KK$.
Then the following regularized adelic formulas for gamma- and beta-functions
are valid (see Eqs (4.18) and (5.4))
$$\prod_{v=1}^\s\G_\infty(\a+i\a_v;\nu_v)\prod_{v=\s+1}^{\s+\u}
\G_{-\infty}(\a+i\a_v;\nu_v)\reg\prod_{v\in F}\G_{q_v}(\a+i\a_v)$$
$$=\kappa\omega (C)[|D|N(J)]^{{1/2}-\a}, \eqno (1.2)$$
$$\prod_{v=1}^\s B_{\infty}(\a+i\a_v,\nu_v;\b+i\b_v,\mu_v)
\prod_{v=\s+1}^{\s+\u}B_{-\infty}(\a+i\a_v,\nu_v;\b+i\b_v,\mu_v)$$
$$\times\reg\prod_{v\in F}B_{q_v}(\a+i\a_v,\b+i\b_v)
={\bf\kappa}\sqrt{|D| N(J)}. \eqno (1.3)$$
where $\G_\infty (\a;\nu)$ and $B_\infty (\a,\nu;\b,\mu)$ are gamma- and
beta-functions of the real numbers field $\QQ_\infty =\RR$ for the characters
(see \Par 2)
$$\t_v(x)={\sgn}^\nu x, \quad \t'_v(x)={\sgn}^\mu x,$$
$$\quad \nu=\nu_v, \mu=\mu_v\in F_2, v=1,2,\ldots,\s;$$
$B_{-\infty}(\a,\nu;\b,\mu)$ is the beta-function of the complex numbers field
$\QQ_{-\infty}=\CC$ for characters (see \Par 2)
$$\t_v(x)=x^\nu (x\bar x)^{-{\nu/2}}, \quad \t'_v(x)
=x^\mu (x\bar x)^{-{\mu/2}},$$
$$\nu=\nu_v, \mu=\mu_v\in Z, \quad v=\s+1,\ldots,\s+\u;$$
and
$$B_{q_v}(\a,\b)=\G_{q_v}(\a)\G_{q_v}(\b)\G_{q_v}(\g), \quad \a+\b+\g=1
 \eqno (1.4)$$
is the beta-function of a local p-field $\KK_v$ with module $q_v$ for 
unramified characters, $v\in F$ (i.e. $v=\s+\u+1,\ldots,$ $\t_v(x)\equiv 1$). 
Here
$$\G_{q_v}(\a)={{1-q_v^{\a-1}}\over{1-q_v^{-\a}}},
 \quad \G_{q_v}(\a)\G_{q_v}(1-\a)=1 \eqno (1.5)$$
is the reduced gamma-function of the field $\KK_v$ (see \Par 2).

For further notations in (1.2) and (1.3) see \Par4 and \Par 5.

For unramified characters formulas (4.19) and (5.5) (special cases of (1.2) and
(1.3)) were obtained in~[9],[10]. Their concrete realizations have been made 
in the following cases: for rational numbers field $\KK=\QQ$ -- in~[11], for 
quadratic fields $\KK=\QQ(\sqrt d)$ -- in~[12], for m-circular fields
$\KK=\QQ(\exp({2\pi i}/m))$ -- in~[10],[9], for some cubic fields
$\KK=\QQ(\root 3\of d)$ -- in~[10].

Adelic formula (1.3) is applied to the four-point Veneziano and
Virasoro-Shapiro amplitudes and their generalizations which connect them with
string (open or closed resp.) amplitudes (see~[1]--[12]).

In this paper we use notations of~[10]. In particular,
$$\delta (|x|-q^\g)=\cases{1, \quad |x|=q^\g \cr 0, \quad |x|\neq q^\g\cr},
 \quad \Omega (q^{-\g}|x|)=\cases{1, \quad |x|\leq q^\g \cr 0,
  \quad |x|>q^\g\cr}.$$
$Z$ is the ring of rational integers; $Z_+$ is the set of natural numbers.

Necessary informations on number theory and analysis can be found in
books~[13]--[20].

The results of this paper were published shortly in Doklady RAN~[21].

This work is supported in part by RFFI Grant 96-01-01008.
\bigskip
\centerline {\BF \Par 2. Local fields of zero characteristic}
\medskip

Such fields are well-known: they are $\RR, \CC$ and algebraic extentions of 
$\QQ_p$ ($p$-fields~[13]). We denote them by $\KK$.

Let $x$ be an element of $\KK$, $|x|$ be a valuation on $\KK$, $dx$ be the
normed Haar mesure on $\KK$. For $p$-field $\KK$ let $q$ be its module, $\ZZ$ 
be the ring of integer elements, $\ZZ^\times$ be the multiplicative subgroup 
of $\ZZ$, $I$ be its maximal ideal.

Let $\chi (x)$ be an additive character of $\KK$, $r$ be its {\it rank} 
(order)~[13],[18]. (For $p$-field $\KK$ the rank $r$ is max of $\g\in Z$ for 
which $\chi (x)\equiv 1, |x|\leq q^\g$.) Therefore the mesure $q^{-{r/2}}dx$ is
self-dual with respect to this character. 

Let
$$\omega (x)=\t (x)|x|^{i\a}, \quad \a\in \RR \eqno (2.1)$$
be a multiplicative character of $\KK$, $\rho$ be its {\it rank}~[13],[18]. 
(For $p$-field $\KK$ the rank $\rho$ is min of $\g\in Z$ for which 
$\omega (x)\equiv 1, x\in\ZZ, 1-x\in I^\g$). 

Now we introduce the following definitions~[18].

{\bf Gamma-function $\G(\a;\t)$} of field $\KK$ for a character $\t$ is the
analytic continuation from domain $\Re\a>0$ of integral
$$\G(\a;\t)=q^{-{r/2}}\int_{\KK}\t (x)|x|^{\a-1}\chi (x)dx. \eqno (2.2)$$

{\bf Beta-function $B(\a,\t;\b,\t')$} of field $\KK$ for characters $\t$ and
$\t'$ is the analytic continuation from tube domain
$\Re\a>0, \Re\b>0, \Re (\a+\b)<1$ of integral
$$B(\a,\t;\b,\t')=\t (-1)\t'(-1)\int_{\KK}\t (x)|x|^{\a-1}\t' (1-x)
|1-x|^{\b-1}dx. \eqno (2.3)$$

{\bf The Mellin transform $\Phi(\a;\t)$} of a function 
$\varphi\in\script{S}(\KK)$ with respect to a character $\t$ of field $\KK$ 
is the analytic continuation from domain $\Re\a>0$ of integral
$$\Phi (\a;\t)=\int_{\KK^\times}\varphi (x)\t (x)|x|^\a d^\times x \eqno(2.4)$$
where $d^\times x$ is the normed Haar mesure on $\KK^\times$.

The following formulas are valid:
$$B(\a,\omega;\b,\omega')=B(\a+i\a',\t;\b+i\b',\t')$$
if $\omega'(x)=\t'(x)|x|^{i\a'}$;
$$B(\a,\t;\b,\t')=q^{r/2}\G(\a;\t)\G(\b;\t')\G(\g;\t'') \eqno (2.5)$$
if $\a+\b+\g=1, \t\t'\t''=1$;
$$\G(\a;\t)\G(1-\a;\bar\t)=\t (-1), \quad \G(\a;\omega)=\G(\a+i\a';\t);
 \eqno(2.6)$$
$$F[\t (x)|x|^{\a-1}](\xi )=\G (\a;\t)\bar\t (\xi )|\xi|^{-\a} \eqno (2.7)$$
where $F[f]=\tilde f$ is the Fourier transform of $f\in\script{S}'(\KK)$ with 
respect to the character $\chi (x)$ and the self-dual mesure $q^{-{r/2}}dx$;
$$\Phi (\a;\t)=\G (\a;\t)\tilde{\Phi} (1-\a;\bar\t)\t(-1),
 \quad \Phi (\a;\omega)=\Phi (\a+i\a';\t) \eqno (2.8)$$
where $\tilde\Phi$ is the Mellin transform of $\tilde\varphi$ with respect to
the character $\bar\t$.

In special local fields the quantities just introduced have the following
particular form.

{\bf Field $\KK=\QQ_\infty=\RR$}. Here $dx$ is the Lebesgue mesure,
$d^\times x=|x|^{-1}dx$, $q=1$, $\chi (x)=\exp (-2\pi ix)$, $r=0$,
$$\t_\nu(x)={\sgn}^\nu x, \quad \nu\in F_2, \quad \rho =0.$$
($F_2$ means the field of residues of $\mod 2$, i.e.~it consists of elements
0 and 1.)
$$\G_\infty (\a;\t_\nu )=\G_\infty (\a;\nu), \quad \nu\in F_2, \eqno (2.9)$$
where
$$\G_\infty (\a;\nu )=i^{-\nu}\pi^{{1/2}-\a}
{\G({{\a+\nu}\over 2})\over\G({{1-\a+\nu}\over 2})},
 \quad \G_\infty (\a;\nu)\G_\infty (1-\a;\nu)=(-1)^\nu; \eqno (2.10)$$
and
$$B_\infty (\a,\t_\nu;\b,{\t'}_\mu)=B_\infty(\a,\nu;\b,\mu),
 \quad \nu, \mu\in F_2, \eqno (2.11)$$
where
$$B_\infty (\a,\nu;\b,\mu)=\G_\infty (\a;\nu)\G_\infty (\b;\mu)
\G_\infty (1-\a-\b;-\nu-\mu). \eqno (2.12)$$

Owing to (2.12) $B_\infty$ is symmetric under permutations of points
$(\a,\nu)$, $(\b,\mu)$, $(\g,\eta)$ on the manifold
$$\a+\b+\g=1, (\a,\b,\g )\in\CC^3, \quad \nu+\mu+\eta=0,
 (\nu,\mu,\eta )\in F_2^3. \eqno (2.13)$$

{\bf Field $\KK=\QQ_{-\infty}=\CC$.} Here $|x|=x\bar x$, 
$dx=|\partial x\wedge\bar\partial x|$, $d^\times x=(x\bar x)^{-1}dx$, 
$q=1$, $\chi (x)=\exp (-2\pi i(x+\bar x))$, $r=0$,
$$\t_\nu (x)=x^\nu (x\bar x)^{-{\nu/2}}, \quad \nu\in Z, \quad \rho =0,$$
$$\G_{-\infty}(\a;\t_\nu )=\G_{-\infty}(\a;\bar{\t}_\nu )
=\G_{-\infty}(\a;\nu )=\G_{-\infty}(\a;-\nu ) \eqno (2.14)$$
where
$$\G_{-\infty}(\a;\nu)=i^{-|\nu |}(2\pi)^{1-2\a}
{{\G(\a+{|\nu|/2})}\over {\G(1-\a+{|\nu|/2})}},$$
$$\G_{-\infty}(\a;\nu)\G_{-\infty}(1-\a;\nu)=(-1)^\nu; \eqno (2.15)$$
and
$$B_{-\infty}(\a,\t_\nu;\b,\t_\mu)=B_{-\infty}(\a,\nu;\b,\mu) \eqno (2.16)$$
where
$$B_{-\infty} (\a,\nu;\b,\mu)=\G_{-\infty} (\a;\nu)\G_{-\infty} (\b;\mu)
\G_{-\infty} (1-\a-\b;-\nu-\mu). \eqno (2.17)$$

Owing to (2.17) $B_{-\infty}$ is symmetric under permutations of points
$(\a,\nu)$, $(\b,\mu)$, $(\g,\eta)$ on the manifold
$$\a+\b+\g=1, (\a,\b,\g )\in\CC^3, \quad \nu+\mu+\eta=0,
 (\nu,\mu,\eta )\in Z^3. \eqno (2.18)$$

{\bf Fields $\KK=\QQ_p(\epsilon), p=2,3,5,\ldots$}. Here $dx$ is the Haar
mesure on $\KK$ normed by condition
$$\int_{|x|\leq 1}dx=1;$$
$$d^\times x=(1-q^{-1})^{-1}|x|^{-1}dx, \quad \int_{|x|=1}d^{\times}x=1$$
where $q=p^f, f\in Z_+$ is the module of $\KK$;
$$\chi (x)=\exp (2\pi i\{x\}_p)$$
where $\{x\}_p$ is the rational part of $x\in\KK$. A multiplicative character
$\t (x)$ on $\KK$ we normalize by the condition
$$\t (\pi)=1 \eqno (2.19)$$
where $\pi$ is the {\it generating element} of $\KK$.

For unramified characters $(\t\equiv 1, \rho=0)$ we have (cf.~[10])
$$\G(\a;1)=q^{-{r/2}}\int_{|x|\leq q^{r+1}}|x|^{\a-1}\chi (x)dx
=q^{(\a-{1/2})r}\G_q(\a) \eqno (2.20)$$
where $\G_q(\a)$ is the reduced gamma-function of field $\QQ_p(\epsilon)$ with
module $q$ defined by the formula (1.5).

(The Eq.~(2.20) follows directly from the $p$-adic integral
$$\int_{|x|\leq q^\g}\chi (\xi x)dx=q^\g\Omega (q^{\g-r}|\xi |) \eqno(2.21).$$

Indeed, it is trivial for $|\xi |\leq q^{r-\g}$. Let $\xi$ be such that 
$|\xi |>q^{r-\g}$. There exists $x_0, |x_0|=q^{r+1}$ such that 
$\chi (\xi x_0)\neq 0$. Then 
$$\int_{|x|\leq q^{\g}}\chi (\xi x)dx=\int_{|x+x_0|\leq q^{\g}}
\chi (\xi (x+x_0))dx=\chi (\xi x_0)\int_{|x|\leq q^{\g}}\chi (\xi x)dx,$$
and (2.21) follows.)

From (2.20) it follows that $\G(\a;1)$ is the Mellin transform of function~[10]
$$(1-q^{-1})q^{-{r/2}}\chi (x)\Omega (q^{-r-1}|x|)\in\script{S}(\KK).
\eqno(2.22)$$
Besides,
$$B(\a,1;\b,1)=B_q(\a,\b)=\G_q(\a)\G_q(\b)\G_q(1-\a-\b). \eqno (2.23)$$

For ramified characters $(\t\not\equiv 0, \rho\geq 1)$ we have~[13]
$$\G(\a;\t)=q^{-{r/2}}\int_{|x|=q^{r+\rho}}\t (x)|x|^{\a-1}\chi (x)dx
=\kappa (\t)q^{(\a-{1/2})(r+\rho)} \eqno (2.24)$$
where
$$\kappa (\t)=q^{\rho/2}\int_{|x|=1}\t(x)\chi (\pi^{-r-\rho}x)dx,
 \quad |\kappa|=1. \eqno (2.25)$$

From (2.24) it follows that $\G(\a;\t)$ is the Mellin transform of function
$$(1-q^{-1})q^{-{r/2}}\chi (x)\delta (|x|-q^{r+\rho})\in\script{S}(\KK).
 \eqno (2.26)$$
\bigskip  
\centerline {\BF \Par 3. Analysis on the adele ring}
\medskip

Let $\KK=\QQ(\epsilon )$ be an algebraic numbers field of degree $n=\s+2\u$,
and $\ZZ$ is its {\it ring of enteger elements} (the maximal order). We denote
$$(\s_v,\KK_v, \quad v=1,2,\ldots)$$
places of field $\KK$~[13]. Here $\KK_v$ is $\RR$ for infinite real places
$v=1,2,\ldots,\s$; $\KK_v$ is $\CC$ for infinite imaginary places
$v=\s+1,\ldots,\s+\u$; $\KK_v$ is p-field $\QQ_p(\s_v(\epsilon ))$ for finite
plces $v=\s+\u+1,\ldots$.

We equip quantities of \Par 2 related to field $\KK_v$ with index $v$, namely:
$x_v$, $|x_v|_v$, $d_vx_v$, $d_v^\times x_v$, $\chi_v(x_v)$, $r_v$. Let
$$\omega_v(x_v)=\t_v(x_v)|x_v|_v^{i\a_v}, \quad \a_v\in\RR \eqno (3.1)$$
be a multiplicative character on $\KK_v$, and $\rho_v$ is its rank. For
p-field $\KK_v (v=\s+\u+1,\ldots)$ let $q_v=p^{f_v}, f_v\in Z_+$  be its
module, $\ZZ_v$ be its ring of enteger elements, $\ZZ_v^\times$ be 
multiplicative subgroup, $I_v$ be maximal ideal of ring $\ZZ_v$, and 
$\Pi_v=\ZZ\cap I_v$ be a prime ideal of the ring $\ZZ$.

Let $A$ be the ring of adeles,
$$X=(x_1,x_2,\ldots,x_v,\ldots), x_v\in\KK_v, x_v\in \ZZ_v, v\geq V>\s+\u),$$
and  $A^\times$ be the group of ideles,
$$X=(x_1,x_2,\ldots,x_v,\ldots, x_v\in\KK_v^\times, x_v\in \ZZ_v^\times,
 v\geq V>\s+\u)$$
of the field $\KK$.

Any character $\chi (X)$ of ring $A$ is represented in the form
$$\chi (X)=\prod_v\chi_v(x_v) \eqno(3.2)$$
where $\chi_v(x_v)$ is the restriction of $\chi (X)$ on subring of adeles of
the form $(0,0,\ldots,x_v,0,\ldots)$. 

We suppose that the character $\chi (X)$ is trivial on $\KK$ considered as
subring of $A$ of principal adeles
$${\bf x}=(\s_v(x), v=1,2,\ldots, x\in\KK ),$$
i.e. $\chi (X)$ is a character of $A/\KK$. As an example of such characters is
the character
$$\chi (X)=\prod_{v=1}^\infty\chi_p ({\Sp}_vx_v).$$

Any character $\omega (X)$ of the group $A^\times$ is represented in the form
$$\omega (X)=\prod_v\omega_v(x_v)=\t (X)\prod_v|x_v|_v^{i\a_v},
 \quad \t (X)=\prod_v\t_v(x_v) \eqno (3.3)$$
where $\omega_v(x_v)$ is defined by Eq.~(3.1). It is the restriction of
$\omega (X)$ on subgroup of ideles of the form $(1,1,\ldots,x_v,1,\ldots)$.

We suppose that the character $\omega (X)$ is trivial on $\KK^\times$ of 
principal ideles, i.e. $\omega (X)$ is a character of 
${{A^\times}/{\KK^\times}}$.

There are necessary and sufficient conditions that a character $\omega (X)$ of
$A^\times$ is trivial on $\KK^\times$~[13]. For instance, for field $\KK=\QQ$
these conditions can be expressed in the explicit form~[18]: {\it a character
$$\omega (X)=\t (X)|x|^{i\a}\prod_{p=2}^\infty |x_p|_p^{i\a_p}, \eqno (3.4)$$
where
$$\t(X)={\sgn}^{\nu}x\prod_{p=2}^\infty\t_p(x_p), \quad \t_p(p)=1,\eqno(3.5)$$
is trivial on $\QQ^\times$ iff
$$\t(-1)=1, \quad p^{i\a_p}=p^{i\a}\t(p). \eqno (3.6)$$}

We recall Eqs~[13],[14]
$$\prod_{v>\s+\u} q_v^{r_v}=|D|, \quad N(J)=\prod_{v\in R}q_v^{\rho_v}.
 \eqno (3.7)$$
where $D$ is discriminant of field $\KK$ and $N(J)$ is the norm of the {\it
principal ideal} 
$$J=\prod_{v\in R}\Pi_v^{\rho_v}$$
of the character $\omega$; $R$ is the set of finite ramified places of $\KK$. 
The Haar mesure
$$\prod_{v=1}^\infty q_v^{-{{r_v}/2}}d_v x_v=|D|^{-{1/2}}dX \eqno (3.8)$$
on $A$ is self-dual with respect to the additive character (3.2) of ring $A$.

Let $V>\s+\u$, and 
$$\varphi (X)=\prod_{v\in F_V}\varphi_v(x_v)\prod_{v\in R}\varphi_v(x_v)
\prod_{v\geq V, v\in F}\Omega (|x_v|_v), \quad \varphi_v\in\script{S}(\KK_v) 
\eqno(3.9)$$
be a standard Schwartz-Bruchat function on $A$ that is 
$\varphi\in\script{S}(A)$~[13]. Here $F$ is the set of finite unramified 
places of field $\KK$, and $F_V=[v\in F: v<V]$.

{\it The Mellin transform $\Phi (\a;\omega)$} of function 
$\varphi\in\script{S}(A)$ with respect to a character $\omega$ of group 
${{A^\times}/{\KK^\times}}$ is the analytic continuation from domain $\Re\a>1$
of the integral 
$$\Phi(\a;\omega)=\int_{A^\times}\varphi (X)\omega (X)|X|^{\a}d^{\times}X
 \eqno(3.10)$$
where $d^\times X$ is the Haar mesure on $A^\times$ defined by the local Haar
mesures $d_v^\times x_v$ of fields $\KK_v$ (\Par 2),
$$d^\times X=\prod_vd_v^\times x_v=\prod_{v=1}^{\s+\u}|x_v|_v^{-1}d_vx_v
\prod_{v>\s+\u}(1-q_v^{-1})^{-1}|x_v|_v^{-1}d_vx_v.$$

The following theorem is valid~[13].

{\BF Theorem} (Tate). {\it For unramified character $\omega$ function
$\Phi (\a;\omega)$ is holomorphic on $\a$ everywhere except simple poles $\a=0$
and $\a=1$ with residues $C\Phi (0)$ and $C'\tilde\Phi (0)$ resp. where numbers
$C$ and $C'$ depend only on field $\KK$. For ramified character $\omega$
function $\Phi (\a;\omega)$ is entire. In addition the folowing relation is
valid,
$$\Phi (\a;\omega)=\tilde{\Phi}(1-\a;\bar\omega) \eqno (3.11)$$
where $\tilde\Phi$ is the Mellin transform of the Fourier transform
$\tilde{\varphi}(\Xi)$ of function $\varphi (X)$ with respect to character
$\bar\omega=\omega^{-1}$.}

In detail the Tate formula (3.11) for function $\varphi$ of form (3.9) can be
rewritten as follows
$$\prod_{v\in F_V}\Phi_v (\a;\omega_v)\prod_{v\in R}\Phi_v (\a;\omega_v)
L_V(\a;\omega)=$$
$$\prod_{v\in F_V}\tilde{\Phi}_v(1-\a;\bar\omega_v)\prod_{v\in R}
\tilde{\Phi}_v(1-\a;\bar\omega_v)L_V(1-\a;\bar\omega). \eqno(3.12)$$
In (3.12) $L_V$ is the $V$-truncated Dirichlet $L$-function,
$$L_V(\a;\omega)=\prod_{v\geq V, v\in F}[1-\lambda (v)q_v^{-\a}]^{-1},
 \lambda (v)=q_v^{-i\a_v}=\omega_v(\pi_v). \eqno (3.13)$$

Recall that the Dirichlet $L$-function $L(\a;\omega)$ is defined by Eq.
$$L(\a;\omega)=\prod_{v\in F}[1-\lambda (v)q_v^{-\a}]^{-1}, \quad \Re>1.
 \eqno (3.14)$$
\bigskip
\centerline {\BF \Par 4. Regularized adelic formulas for gamma-functions}
\medskip

Our goals are to describe a regularization and to sketch a proof of formula
(1.2) for gamma-functions. To this end, we specify the Tate formula (3.12) by 
a special choice of standard function $\varphi (X)$ (see (3.9)), similar to 
what has been done in~[9]--[12].

For $v=1,2,\ldots,\s$ we put
$$\varphi_v (x)=x^\nu\exp (-\pi x^2),\quad \tilde{\varphi}_v (\xi )
=i^{-\nu}\varphi_v (\xi), \quad \nu=\nu_v\in F_2. \eqno (4.1)$$
Then
$$\Phi_v(\a;\omega_\nu)=G_{\infty}(\a+i\a_v+\nu), \quad 
\tilde{\Phi}_v(\a;\bar{\omega}_\nu)=i^{-\nu}G_\infty(\a-i\a_v+\nu),
 \eqno (4.2)$$
where~[13]
$$G_\infty (\a)=\int_{-\infty}^\infty |x|^{\a-1}\exp (-\pi x^2)dx
=\pi^{-{\a/2}}\G({\a/2}).$$
Note owing to (2.10) that
$$\G_\infty (\a;\nu)=i^{-\nu}{{G_\infty (\a+\nu)}\over{G_\infty (1-\a+\nu)}},
 \quad \nu\in F_2. \eqno (4.3)$$

For $v=\s+1,\ldots,\s+\u$ we put
$$\varphi_v (x)=x^{\nu}\exp (-2\pi x\bar x), \nu\geq 0, \quad \varphi_v (x)
=\bar {x}^{|\nu|}\exp (-2\pi x\bar x), \nu<0,$$
$$\tilde{\varphi}_v(\xi )=i^{-\nu}\bar{\varphi}_v(\xi), \quad \nu=\nu_v\in Z.$$
Then
$$\Phi_v(\a;\omega_\nu)=G_{-\infty}(\a+i\a_v+{{|\nu|}/2}),$$
$$\tilde{\Phi}_v(\a;\bar{\omega}_\nu)
=i^{-\nu}G_{-\infty}(\a-i\a_v+{{|\nu|}/2}) \eqno (4.4)$$
where~[13]
$$ G_{-\infty}(\a)=\int_{\CC} (x\bar x)^{\a-1}\exp (-2\pi x\bar x)dx
=(2\pi)^{1-\a}\G(\a).$$
Note owing to (2.15) that
$$\G_{-\infty}(\a;\nu)
=i^{-|\nu |)}{{G_{-\infty}(\a+{|\nu|/2})}\over{G_{-\infty}(1-\a+{|\nu|/2})}},
 \quad \nu\in Z. \eqno (4.5)$$

For $v\in F_V$ we take function $\varphi_v$ in form (2.22) (for field 
$\KK_v$). Then owing to (2.20) we have
$$\Phi_v(\a;\omega_v)=\G_v(\a+i\a_v;1)
=q_v^{(\a+i\a_v-{1/2})r_v}\G_{q_v}(\a+i\a_v) \eqno (4.6)$$
and owing to (2.8) and (4.6)
$$\tilde{\Phi}_v(\a;\bar{\omega}_v)=1. \eqno (4.7)$$

For $v\in R$ we take function $\varphi_v$ in form (2.26) (for field $\KK_v$).
Then owing to (2.24) we have
$$\Phi_v(\a;\omega_v)=\G_v(\a+i\a_v;\t_v)
=\kappa_vq_v^{(\a+i\a_v-{1/2})(r_v+\rho_v)}, \eqno (4.8)$$
where a number $\kappa_v=\kappa (\t_v)$ is defined by Eq. (2.25),
$$\kappa_v(\t_v)=q_v^{{\rho_v}/2}\int_{|x|_v=1}\t_v (x)
\chi (\pi_v^{-r_v-\rho_v}x)d_vx, \quad |\kappa_v|=1. \eqno (4.9)$$
At last, owing to (2.8) and (4.8) we have
$$\tilde{\Phi} (\a;\bar{\omega}_v)=\t_v(-1). \eqno (4.10)$$

Substituting expressions (4.2),(4.4),(4.6)-(4.8) and (4.10) in the Tate formula
(3.12) we get for all $V>\s+\u$ the following Eq.
$$\prod_{v=1}^{\s}G_\infty (\a+i\a_v+\nu_v)\prod_{v=\s+1}^{\s+\u}
G_{-\infty} (\a+i\a_v+{{|\nu_v|}/2})\prod_{v\in F_V}q_v^{(\a+i\a_v-{1/2})r_v}
\G_{q_v}(\a+i\a_v)$$
$$\times\prod_{v\in R}\kappa (\t_v)q_v^{(\a+i\a_v-{1/2})(r_v+\rho_v)}
L_V(\a;\omega)=\prod_{v=1}^\s i^{-\nu_v}G_\infty(1-\a-i\a_v+\nu_v)$$
$$\times\prod_{v=\s+1}^{\s+\u}i^{-\nu_v}G_{-\infty}(1-\a-\a_v+{{|\nu_v|}/2})
\prod_{v\in R}\t_v(-1)L_V(1-\a;\bar\omega). \eqno (4.11)$$

Now we divide both sides of Eq.~(4.11) on
$$\prod_{v=1}^{\s}G_\infty (1-\a-i\a_v+\nu_v)
\prod_{v=\s+1}^{\s+\u}G_{-\infty}(1-\a-i\a_v+{{|\nu_v|}/2}).$$
As a result, using Eqs (4.3),(4.5) and (3.7) (cf.~[13],[10]), we get the
following key formula for gamma-functions
$$\prod_{v=1}^\s\G_\infty (\a+i\a_v;\nu_v)\prod_{v=\s+1}^{\s+\u}
\G_{-\infty}(\a+i\a_v;\nu_v)\prod_{v\in F_V}\G_{q_v}(\a+i\a_v)L_V(\a;\omega)$$
$$=\kappa\omega (C)[|D|N(J)]^{{1/2}-\a}L_V(1-\a;\bar\omega) \eqno (4.12)$$
for all $V>\s+\u$. In (4.12)
$$\kappa=\prod_{v=1}^\infty\kappa_v, \quad |\kappa|=1, \eqno (4.13)$$
where $\kappa_v=i^{-\nu_v}, v=1,2,\ldots,\s,$ $\kappa_v=i^{-|\nu_v |},$ 
$v=\s+1,\ldots ,\s+\u, \quad \kappa_v=1, v\in F$,
$\kappa_v=\t_v(-1)\bar{\kappa}(\t_v), v\in R$ (see (4.9)); and also
$$C=\prod_{v>\s+\u} \pi_v^{r_v+\rho_v},
 \quad \omega (C)=\prod_{v>\s+\u} q_v^{-i\a_v(r_v+\rho_v)}. \eqno (4.14)$$

As Euler's function $\G (\a)$ does not vanish everywhere, so both sides of 
Eq.~(4.12) are entire functions provided the set $R$ is not empty.

Formula (4.12) is valid also for $V=\s+\u+1$ ($F_V=F$). In this case owing to
(3.14) it gives the functional Eq.~for the Dirichlet $L$-function~[13],
$$\prod_{v=1}^\s\G_\infty (\a+i\a_v;\nu_v)\prod_{v=\s+1}^{\s+\u}
\G_{-\infty} (\a+i\a_v;\nu_v)L(\a;\omega)$$
$$=\kappa\omega (C)[|D|N(J)]^{{1/2}-\a}L(1-\a;\bar\omega). \eqno(4.15)$$

By virtue of (3.13) for $\Re\a<0$
$$\lim_{V\to\infty}L_V(1-\a;\bar\omega)=1, \eqno (4.16)$$
and thus by (4.12) there exists
$$\lim_{V\to\infty}\prod_{v\in F_V}\G_{q_v}(\a+i\a_v)L_V(\a;\omega)$$
which we will denote by
$$\reg\prod_{v\in F}\G_{q_v}(\a+i\a_v), \eqno (4.17)$$,
and this limit satisfies Eq.
$$\prod_{v=1}^\s\G_\infty(\a+i\a_v;\nu_v)\prod_{v=\s+1}^{\s+\u}
\G_{-\infty}(\a+i\a_v;\nu_v)\reg\prod_{v\in F}\G_{q_v}(\a+i\a_v)$$
$$=\kappa\omega (C)[|D|N(J)]^{{1/2}-\a}. \eqno (4.18)$$

For $\Re\geq 0$ the regularization expression (4.17) is defined by Eq.~(4.18) 
as the analytic continuation on $\a$ from domain $\Re\a<0$ (as a meromorphic 
function).

{\it Eq.~(4.18) is regularized adelic formula for gamma-functions of all 
completions of field} $\KK$.

For unramified characters $(\nu_v=0, \a_v=0, \kappa=1, \omega (C)=1, N(J)=1)$
Eq.~(4.18) takes the form~[9],[10]
$$\G_\infty^\s (\a)\G_{-\infty}^\u (\a)
\reg\prod_{p=2}^\infty\prod_{j=1}^{m_p}\G_{q_{pj}}(\a)=|D|^{{1/2}-\a}
 \eqno(4.19)$$
where $m_p$ is number of prime divisors of field $\KK$ entering in
decomposition of prime number $p$.

For field $\KK=\QQ \quad (\s=1, \u=0, \a_\infty=0, r_v=0, D=1)$ Eq.~(4.18)
takes the form
$$\G_\infty (\a;\nu)\reg\prod_{p\in F}\G_p(\a+i\a_p)
=\kappa\omega (C)N^{{1/2}-\a}(J), \eqno (4.20)$$
where numbers $\nu, \a_p, \kappa, \omega (C), N(J)$ are defined by Eqs (3.6),
(3.7), (4.13), (4.14). Note that owing to (3.6) and (1.5) in this case
$$\G_p(\a+i\a_p)={{1-\t(p)p^{\a-1})}\over{1-\bar\t(p)p^{-\a})}},
 \quad p=2,3,5,\ldots. \eqno (4.21)$$
\bigskip
\centerline{\BF \Par 5. Regularized adelic formulas for beta-functions}
\medskip

Let characters $\omega$ and $\omega'$ be defined on
${{A^\times}/{\KK^\times}}$. Suppose that ranks $\rho_v, \rho'_v$ and
$\rho''_v$ of characters $\omega_v,$ $\omega'_v$ and
$\omega''_v=\omega_v\omega'_v$ (of form (3.1)) of fields $\KK_v$ are equal,
$\rho_v=\rho'_v=\rho''_v$, for all $v\in R$. In this case the quantities
$N(J)$ and $C$ and the sets $F$ and $R$ turn out to be the same for 
the characters $\omega$, $\omega'$ and $\omega''$.

Our goal is to prove the regularized formula (1.3) for beta-functions.

If we multiply Eq.~(4.18) on the same Eq.~in which magnitudes $\a, \a_v,
\nu_v, \kappa_v$ and $\omega$ are changed by $\b, \b_v, \mu_v, \kappa'_v$ and
$\omega'$ resp., and the Eq.~obtained we divide on Eq.~(4.18) in which just
mentioned magnitudes are changed by $\a+\b, \a_v+\b_v, \nu_v+\mu_v,
\kappa''_v$ and $\omega''$ resp. As a result by using formulas (2.10),
(2.12), (2.15), (2.17), (2.23) and (1.5) we get for $V>\s+\u$ the
following Eq.
$$\prod_{v=1}^\s B_\infty(\a+i\a_v,\nu_v;\b+i\b_v,\mu_v)
\prod_{v=\s+1}^{\s+\u}B_{-\infty}(\a+i\a_v,\nu_v;\b+i\b_v,\mu_v)$$
$$\times\prod_{v\in F_V}B_{q_v}(\a+i\a_v,\b+i\b_v)L_V(\a;\omega)
L_V(\b;\omega')L_V^{-1}(\a+\b;\omega'')$$
$$={\bf\kappa}\sqrt{|D| N(J)}L_V(1-\a;\bar\omega)L_V(1-\b;{\bar\omega}')
L_V^{-1}(1-\a-\b;{\bar\omega}''), \eqno (5.1)$$
where
$${\bf\kappa}=\prod_{v\in R}\kappa_v\kappa'_v\bar{\kappa}''_v,
 \quad |{\bf\kappa}|=1. \eqno (5.2)$$

Passing on to limit as $V\to\infty$ in Eq.~(5.1) for $\Re\a<0$ and $\Re\b<0$,
using (4.16) and denoting
$$\lim_{V\to\infty}\prod_{v\in F_V}B_{q_v}(\a+i\a_v,\b+i\b_v)
L_V(\a;\omega)L_V(\b;\omega')L_V^{-1}(\a+\b;\omega'')$$
$$=\reg\prod_{v\in F}B_{q_v}(\a+i\a_v,\b+i\b_v), \eqno (5.3)$$
we get Eq.
$$\prod_{v=1}^\s B_{\infty}(\a+i\a_v,\nu_v;\b+i\b_v,\mu_v)
\prod_{v=\s+1}^{\s+\u}B_{-\infty}(\a+i\a_v,\nu_v;\b+i\b_v,\mu_v)$$
$$\times\reg\prod_{v\in F}B_{q_v}(\a+i\a_v,\b+i\b_v)
={\bf\kappa}\sqrt{|D| N(J)}. \eqno (5.4)$$

So we have just obtained the following

{\BF Theorem.} {\it Let $\KK$ be an algebraic numbers field of degree
$n=\s+2\u$, $\omega, \omega'$ and $\omega''=\omega\omega'$ be characters on
${{A^\times}/{\KK^\times}}$ and their local ranks $\rho_v$, $\rho'_v$ and
$\rho''_v$ be equal, $\rho_v={\rho}'_v=\rho''_v, v\in R$. Then for
beta-functions of all completions of field $\KK$ the regularized adelic formula
(5.4) is valid.}

For unramified characters formula (5.4) takes the form~[1],[2]
$$B_\infty^\s (\a,\b)B_{-\infty}^\u (\a,\b)
\reg\prod_{p=2}^\infty\prod_{j=1}^{m_p}B_{q_{pj}}(\a,\b)=\sqrt {|D|}.
 \eqno (5.5)$$

For field $\KK=\QQ$ formula (5.4) takes the following explicit form
$$B_\infty(\a,\nu;\b,\mu)\reg\prod_{p\in F}B_p(\a+i\a_p,\b+i\b_p)
={\bf \kappa}\sqrt {N(J)} \eqno (5.6)$$
where owing to (3.6) and (3.5)
$$\t(-1)=\t'(-1)=1, \quad p^{i\a_p}=\t(p), \quad p^{i\b_p}
=\t'(p), p\in F. \eqno (5.7)$$

Similar, in a more complicated way, regularized adelic formulas can be written
out for cases when local ranks $\rho_v, \rho'_v, \rho''_v$ of characters 
$\omega, \omega', \omega''$ are different.
\vskip 12 truemm
\centerline{\BF Refrences}
\bigskip

 \item{1}  P.~G.~O.~Freund and E.~Witten, Adelic String Amplitudes, 
           \underbar{Phys. Lett. B}, \underbar{199}:191--194 (1987).

 \item{2.} I.~V.~Volovich, Harmonic Analysis and p-Adic Strings, 
           \underbar{Lett. Math. Phys}, \underbar{16}:61--67 (1988).

 \item{3.} P.~G.~O.~Freund and M.~Olson, Non-Archimedian Strings, 
           \underbar{Phys. Lett. B}, \underbar{199}:186--190 (1987).

 \item{4.} I.~V.~Volovich, p-Adic String, \underbar{Class. Quantum Grav.},
           \underbar{4}:L83--L87 (1987).

 \item{5.} I.~Ya.~Aref'eva, B.~G.~Dragovi\v c and I.~V.~Volovich, On the 
           Adelic String Amplitudes, \underbar{Phys. Lett. B}, 
	   \underbar{209}:445--450 (1988).

 \item{6.} I.~Ya.~Aref'eva, B.~G.~Dragovi\v c and I.~V.~Volovich, Open and
           Closed p-Adic Strings and Quadratic Extentions of Number Fields,
	   \underbar{Phys. Lett. B}, \underbar{212}:283--291 (1988).

 \item{7.} P.~H.~Frampton, Y.~Okada and M.~R.~Ubriaco, On Adelic Formulas for
           p-Adic String, \underbar{Phys. Lett. B}, \underbar{213}:260--262 
	   (1988).

 \item{8.} P.~H.~Frampton, and H.~Nishino, Theory of p-Adic Closed Strings, 
           \underbar{Phys. Rev. Lett.}, \underbar{62}:1960--1963 (1989).

 \item{9.} V.~S.~Vladimirov and T.~M.~Sapuzhak, Adelic Formulas for String
           Amplitudes in Fields of Algebraic Numbers, 
	   \underbar{Lett. Math. Phys.}, \underbar{37(2)}:233--242 (1996).

 \item{10.}V.~S.~Vladimirov, Adelic Formulas for Gamma- and Beta-functios of
           Completions of Algebraic Number Fields and Their Applicatios to
           String Amplitudes, \underbar{Izv. RAN, ser.~math.}, \underbar{60}, 
	   \No1:63--86 (1996)(in Russian).

 \item{11.}V.~S.~Vladimirov, On the Freund-Witten Adelic Formula for Veneziano
           Amplitudes, \underbar{Lett. Math. Phys.}, \underbar{27}:123--131 (1993).

 \item{12.}V.~S.~Vladimirov, The Adelic Freund-Witten Formulas for the Veneziano
           and Virasoro-Shapiro Amplitudes, \underbar{Uspechi Math. Nauk.},
	   \underbar{48}, \No6:3--38 (1993) (in Russian).
           
 \item{13.}A.~Weil. \underbar{Basic Numbers Theory}, Springer, Berlin etc, 
           1967.

 \item{14.}Z.~I.~Borevich and I.~R.~Shafarevich. 
           \underbar{Th\'eory de nombres}, Gauthier-Villars, Paris, 1967.

 \item{15.}S.~Lang. \underbar{Algebraic Numbers}, Adisson-Wesley, London, 1964.

 \item{16.}W.~H.~Schikhof. \underbar{Ultrametric Calculus. An Introduction to p-Adic}
           \underbar{Analysis}, Cambridge Univ. Press, 1984,  XI+306p.

 \item{17.}N.~Koblitz. \underbar{p-Adic Numbers, p-Adic Analysis, and Zeta-functions}, 
           Springer, Berlin etc, 1977.

 \item{18.}I.~M.~Gel'fand, M.~I.~Graev and I.~I.~Pjatetskii-Shapiro.
           \underbar{Representations Theory}
	   \underbar{and Automorphic Functions}, Saunders, Philadelphia, 1969.

 \item{19.}V.~S.~Vladimirov, I.~V.~Volovich and E.~I.~Zelenov. 
           \underbar{p-Adic Analysis and}
	   \underbar{Mathematical Physics}, World Scientific, Singapore, 1994, 
	   XVIII+319p.

 \item{20.}Andrei Khrennikov. \underbar{p-Adic Valued Distributions in Mathematical}
           \underbar{Physics}, Kluver, Dordrecht etc, 1994, XIV+264p.

 \item{21.}V.~S.~Vladimirov. Adelic Formulas for Gamma- and Beta-functions in 
           Fields of Algebraic Numbers, \underbar{Doklady Mathematics}, 
	   \underbar{347(1)}:11--15 (1996).
\bye